\documentclass[numbers,preprint,12pt]{elsarticle}
\journal{}

\usepackage{amsmath,amsfonts,amssymb,amsthm,graphicx,hyperref,url}
\usepackage{mathtools,enumerate,xcolor,appendix,dsfont,geometry}
\usepackage{graphicx}
\newcommand{\meddot}{\mathbin{\vcenter{\hbox{\scalebox{0.5}{$\bullet$}}}}}

\theoremstyle{plain}
\newtheorem{theorem}{Theorem}
\newtheorem{proposition}{Proposition}
\newtheorem{lemma}{Lemma}
\newtheorem{corollary}{Corollary}
\newtheorem{remark}{Remark}

\newcommand{\N}{\mathbb{N}}
\newcommand{\R}{\mathbb{R}}
\newcommand{\EE}{\mathsf{E}} 
\newcommand{\Var}{\mathsf{Var}} 
\newcommand{\bb}[1]{\boldsymbol{#1}}
\newcommand{\rd}{\mathrm{d}}

\allowdisplaybreaks

\begin{document}

\begin{frontmatter}

\title{Results related to the Gaussian product inequality conjecture \\for mixed-sign exponents in arbitrary dimension}%

\author[a1]{Guolie Lan}\ead{langl@gzhu.edu.cn}%
\author[a2]{Fr\'ed\'eric Ouimet}\ead{frederic.ouimet2@uqtr.ca}%
\author[a3]{Wei Sun}\ead{wei.sun@concordia.ca}%

\address[a1]{School of Economics and Statistics, Guangzhou University, China}
\address[a2]{D\'epartement de math\'ematiques et d'informatique, Universit\'e du Qu\'ebec \`a Trois-Rivi\`eres, Canada}
\address[a3]{Department of Mathematics and Statistics, Concordia University, Canada}

\begin{abstract}
This paper studies Gaussian product inequalities (GPIs) for centered Gaussian random vectors when positive and negative exponents appear simultaneously. We prove a quantitative lower bound for arbitrary mixed-sign patterns, conditional on the validity of a lower-dimensional GPI, and an upper bound for products of negative powers with convex functions of the remaining Gaussian components. The latter yields an explicit upper bound for mixed-sign products whenever the total positive-side exponent is at least $1/2$.
\end{abstract}


\end{frontmatter}

\section{Introduction\label{sec:intro}}

The real linear polarization constant conjecture grew out of the work of \citet{MR1636556} on obtaining lower bounds for the norms of products of polynomials on real Banach spaces, and, in the linear case, for products of bounded linear functionals. A decade later, \citet{MR2385646} formulated, in connection with the problem restricted to linear functionals on a real Euclidean space, the following Gaussian product inequality (GPI) conjecture: for every real centered Gaussian random vector $\bb{X} = (X_1, \ldots, X_d)$ with $d \in \N$ and every $m \in \N$,
\begin{equation}\label{eq:GPI.original}
\EE\left[\prod_{i=1}^d X_i^{2 m}\right] \geq \prod_{i=1}^d \EE\left[ X_i^{2 m} \right].
\end{equation}
Inequality~\eqref{eq:GPI.original} implies the real linear polarization constant conjecture for real Hilbert spaces \citep{MR3425898}. It is also closely linked to the $U$-conjecture, which asserts that, if two polynomials in the components of a standard Gaussian random vector are independent, an orthogonal transformation can be found so that each depends on a disjoint subset of variables; see, e.g., \citet{MR0346969} and \citet{MR3425898}.

Building on this formulation, \citet{MR2886380} proposed a strengthened conjecture allowing arbitrary positive exponents. Writing these exponents as $2\bb{\nu} = (2\nu_1, \ldots, 2\nu_d)$ and restricting attention to $\bb{\nu} \in (0,\infty)^d$, the resulting statement, denoted by $\mathbf{GPI}_d(\bb{\nu})$, reads
\begin{equation}\label{eq:GPI}
\mathbf{GPI}_d(\bb{\nu}) : \quad \EE\left[ \prod_{i=1}^d |X_i|^{2 \nu_i} \right] \geq \prod_{i=1}^d \EE\left[ |X_i|^{2 \nu_i} \right].
\end{equation}
The conjecture is settled for $d=2$ because $(|X_1|,|X_2|)$ is multivariate totally positive of order~2 \citep{MR628759}. In dimension~$3$, numerous special cases have been confirmed:
\begin{enumerate}[(a)]\setlength\itemsep{0em}
    \item \citet{MR4052574}: $\mathbf{GPI}_3(p,p,q)$ for all $p,q \in \N$;
    \item \citet{MR4593134}: $\mathbf{GPI}_3(1,p,q)$ for all $p,q \in \N$;
    \item \citet{MR4760098}: $\mathbf{GPI}_3(2,3,p)$ for all $p \in \N$;
    \item \citet{MR4661091}: $\mathbf{GPI}_3(p,q,r)$ for all $p,q,r \in \N$;
    \item \citet{MR4798604}: $\mathbf{GPI}_3(p,q,r)$ for $p \in \N$ and $q,r\in (0,\infty)$;
    \item \citet{MR4794515}: $\mathbf{GPI}_3(p,q,r)$ for all $p,q,r \in (1/2)\N$.
\end{enumerate}
For higher $d$, additional results appear under various covariance constraints; see, for instance, \citet{MR4466643}, \citet{MR4445681}, and \citet{MR4554766}.

Interest has also turned to scenarios where some exponents in the product in~\eqref{eq:GPI} are allowed to be negative (subject to the expectations being finite). When every exponent is negative, \citet{MR3278931} proved the inequality and even extended it to multivariate gamma random vectors (also called permanental vectors). Wei's results and related developments extend to traces and determinants of disjoint principal blocks of Wishart random matrices, as shown by \citet{MR4538422} and \citet{MR4822653,MR4862164}.

In the Gaussian setting, the mixed-sign case, where positive and negative exponents occur in the same product, remains largely unresolved. The two-dimensional situation was settled by \citet{MR4471184} (tighter quantitative bounds have also been proved by \citet{MR4530374}), and partial progress for general $d$ was achieved recently by \citet{MR4666255}. Letting $\Sigma = (\sigma_{ij})_{1\leq i,j \leq d}$ be the covariance matrix of $\bb{X}$, they established
\begin{equation}\label{eq:Zhou.et.al.Thm.1.1}
\EE\left[|X_1|^{-2\nu_1} \prod_{i=2}^d X_i^2\right] \geq \left\{\prod_{i=2}^d \left[1 - \frac{\sigma_{1i}^2}{\sigma_{11} \sigma_{ii}}\right]\right\} \EE \left[ |X_1|^{-2\nu_1} \right] \prod_{i=2}^d \EE[X_i^2],
\end{equation}
for any $\nu_1 \in [0,1/2)$, and
\begin{equation}\label{eq:Zhou.et.al.Thm.1.3}
\EE\left[\left(\prod_{i=1}^{d-1} |X_i|^{-2\nu_i}\right) |X_d|^{2\nu_d}\right] \leq \EE\left[\prod_{i=1}^{d-1} |X_i|^{-2\nu_i}\right] \EE \left[ |X_d|^{2\nu_d} \right],
\end{equation}
for any $\nu_1,\ldots,\nu_{d-1} \in [0,1/2)$ and $\nu_d \in (0,\infty)$.

These two estimates already display the contrasting roles of positive and negative powers in mixed-sign problems. In~\eqref{eq:Zhou.et.al.Thm.1.1}, when one negative power is present, the inequality keeps the lower-bound direction of the base GPI~\eqref{eq:GPI}, but a covariance-dependent correction factor appears; no such factor is present in the conjecture~\eqref{eq:GPI} for all positive exponents. In~\eqref{eq:Zhou.et.al.Thm.1.3}, by contrast, the mixed product is bounded in the opposite direction from the base GPI when negative powers are paired with one positive power.

The present paper investigates the general mixed-sign problem. Theorem~\ref{thm:lower.bound} extends \eqref{eq:Zhou.et.al.Thm.1.1} from one negative-exponent component to an arbitrary negative-exponent block, conditional on the validity of a lower-dimensional GPI for the positive-exponent block. Theorem~\ref{thm:upper.bound.convex} extends the upper-bound direction of~\eqref{eq:Zhou.et.al.Thm.1.3} from a single positive-side component to a convex function of an arbitrary positive-side block, with an arbitrary negative-side block, and yields explicit upper bounds for mixed-sign products when the total positive-side exponent is at least $1/2$. Except for the narrow remaining product case with at least two positive-side components and total positive-side exponent in $(0,1/2)$, discussed in Remark~\ref{rem:2}, the corresponding mixed-sign lower and upper estimates follow from Theorems~\ref{thm:lower.bound} and~\ref{thm:upper.bound.convex}, together with~\eqref{eq:Zhou.et.al.Thm.1.3}, although the lower-bound estimates remain subject to stated lower-dimensional GPI assumption.

The sign split is important in these extensions. Negative powers are treated through Laplace-transform representations (see \eqref{eq:completely.monotonicity}), which turn them into Gaussian exponential tilts. After renormalization, the variables carrying positive powers have a reduced covariance matrix related to the Schur complement of the block indexed by the negative powers. The positive powers are then handled either by a lower-order GPI, for the lower bound, or by a convex-order comparison, for the upper bound. This is why the mixed-sign setting can lead both to correction factors, which measure the loss from marginal to conditional variances, and to inequalities in the direction opposite to the base GPI.

The remainder of the paper is organized as follows. Section~\ref{sec:results} states the main results; the proofs are given in Section~\ref{sec:proofs}. In the interest of self-containment, some technical lemmas used in the proofs are relegated to \ref{app:tech.lemmas}.

\section{Results}\label{sec:results}

The first theorem below generalizes~\eqref{eq:Zhou.et.al.Thm.1.1} from one Gaussian component having a negative exponent to an arbitrary subset of the Gaussian components having negative exponents. The positive-side exponents $\nu_i$ are not restricted to be equal to $1$, but this comes at the cost of assuming the validity of a lower-order GPI. The appearance of $\Sigma/\Sigma_{\mathcal{J}\mathcal{J}}$ reflects the fact that the negative-exponent block is first separated, leaving the positive-exponent block with a covariance matrix bounded below by its conditional covariance. Thus, in the presence of genuinely negative powers, the theorem supports the natural conjectural picture that the lower-bound direction of the base GPI can still be retained after inserting explicit Schur-complement correction factors, provided the required GPI is available on the positive side. However, this also means that if the GPI assumption has already been proved, the conclusion of the theorem holds `unconditionally'. For a comprehensive list of results related to the GPI, refer to Section~1 of \citet{MR4822653}.

\begin{theorem}[Lower bound]\label{thm:lower.bound}
Let $\bb{X}\sim \mathcal{N}_d(\bb{0}_d, \Sigma)$ for some positive definite matrix $\Sigma$, and let $\emptyset\neq \mathcal{J}\subsetneq \{1,\ldots,d\}$ be given. Assume that $\nu_j\in [0,1/2)$ for all $j\in \mathcal{J}$ and $\nu_i\in (0,\infty)$ for all $i\in \mathcal{J}^{\complement}$. Then, conditionally on $\mathbf{GPI}_{d-|\mathcal{J}|}((\nu_i)_{i\in \mathcal{J}^{\complement}})$ being true, we have
\[
\EE\left[\prod_{j\in\mathcal{J}}|X_j|^{-2\nu_j} \times \prod_{i\in\mathcal{J}^{\complement}}|X_i|^{2\nu_i}\right]
\geq \EE\left[\prod_{j\in \mathcal{J}} |X_j|^{-2\nu_j}\right] \times \prod_{i\in\mathcal{J}^{\complement}} \left\{\frac{(\Sigma/\Sigma_{\mathcal{J}\mathcal{J}})_{ii}}{\Sigma_{ii}}\right\}^{\nu_i} \EE\left[|X_i|^{2\nu_i}\right],
\]
where $\mathcal{J}^{\complement} \equiv \{1,\ldots,d\}\setminus \mathcal{J}$, the matrix $\Sigma_{\mathcal{I}\mathcal{I}'}$ denotes $\Sigma$ restricted to the rows and columns indexed by $\mathcal{I}$ and $\mathcal{I}'$, respectively, and $\Sigma/\Sigma_{\mathcal{J}\mathcal{J}} = \Sigma_{\mathcal{J}^{\complement}\mathcal{J}^{\complement}} - \Sigma_{\mathcal{J}^{\complement}\mathcal{J}} \Sigma_{\mathcal{J}\mathcal{J}}^{-1} \Sigma_{\mathcal{J}\mathcal{J}^{\complement}}$ is the Schur complement of $\Sigma_{\mathcal{J}\mathcal{J}}$ in $\Sigma$, with rows and columns indexed by $\mathcal{J}^{\complement}$; in $(\Sigma/\Sigma_{\mathcal{J}\mathcal{J}})_{ii}$, the index $i$ is understood as the original index in $\mathcal{J}^{\complement}$.
\end{theorem}

\begin{remark}\label{rem:1}
In Theorem~\ref{thm:lower.bound}, the factor $(\Sigma/\Sigma_{\mathcal{J}\mathcal{J}})_{ii} / \Sigma_{ii}$ is less than or equal to $1$ because $\Sigma/\Sigma_{\mathcal{J}\mathcal{J}} \preceq \Sigma_{\mathcal{J}^{\complement}\mathcal{J}^{\complement}}$ in the Loewner order. Since, in the Gaussian case, $\Sigma/\Sigma_{\mathcal{J}\mathcal{J}}$ is the conditional covariance matrix of $(X_i)_{i\in \mathcal{J}^{\complement}}$ given $(X_j)_{j\in \mathcal{J}}$, this correction factor measures the variance loss caused by conditioning on the variables in the negative-side block. In the special case where $\mathcal{J} = \{1\}$ and $\mathbf{GPI}_{d-1}(1,\ldots,1)$ is known to hold \citep{MR2385646}, we recover~\eqref{eq:Zhou.et.al.Thm.1.1} with
\[
\frac{(\Sigma/\Sigma_{\mathcal{J}\mathcal{J}})_{ii}}{\Sigma_{ii}} = \frac{\sigma_{ii} - \sigma_{i1} \sigma_{11}^{-1} \sigma_{1i}}{\sigma_{ii}} = 1 - \frac{\sigma_{1i}^2}{\sigma_{11} \sigma_{ii}}, \quad i\in \{2,\ldots,d\}.
\]
\end{remark}

It follows from the work of \citet{MR3278931} that the GPI holds when all the exponents are negative. Combining this fact with Theorem~\ref{thm:lower.bound} gives the following product version of the lower bound, still under the same lower-dimensional GPI assumption.

\begin{corollary}\label{cor:lower.bound}
In the setting of Theorem~\ref{thm:lower.bound}, and still conditionally on $\mathbf{GPI}_{d-|\mathcal{J}|}((\nu_i)_{i\in \mathcal{J}^{\complement}})$ being true, we have
\[
\EE\left[\prod_{j\in\mathcal{J}}|X_j|^{-2\nu_j} \times \prod_{i\in\mathcal{J}^{\complement}}|X_i|^{2\nu_i}\right]
\geq \prod_{j\in \mathcal{J}} \EE\left[|X_j|^{-2\nu_j}\right] \times \prod_{i\in\mathcal{J}^{\complement}} \left\{\frac{(\Sigma/\Sigma_{\mathcal{J}\mathcal{J}})_{ii}}{\Sigma_{ii}}\right\}^{\nu_i} \EE\left[|X_i|^{2\nu_i}\right].
\]
\end{corollary}

Next, in Theorem~\ref{thm:upper.bound.convex} below, we give a convex-order upper bound in the spirit of~\eqref{eq:Zhou.et.al.Thm.1.3}. The positive-side factor is replaced by a convex function of the Gaussian components indexed by $\mathcal{J}^{\complement}$. This replacement, in turn, gives explicit upper bounds for mixed-sign products; see Corollary~\ref{cor:upper.bound}.

The upper-bound result reflects the other side of the same sign split. Once the negative powers have been absorbed into the Gaussian exponential tilt, the remaining covariance is smaller than $\Sigma$ in the Loewner order. For convex functions of the components indexed by $\mathcal{J}^{\complement}$, this covariance reduction gives an upper bound through convex order, without any assumption on the GPI. This is the mechanism behind the opposite direction in~\eqref{eq:Zhou.et.al.Thm.1.3}.

\begin{theorem}[Upper bound]\label{thm:upper.bound.convex}
Let $\bb{X}\sim \mathcal{N}_d(\bb{0}_d, \Sigma)$ for some positive definite matrix $\Sigma$, and let $\emptyset\neq \mathcal{J}\subsetneq \{1,\ldots,d\}$ be given. Assume that $\nu_j\in [0,1/2)$ for all $j\in \mathcal{J}$. Then, for any convex function $\psi:\R^{d-|\mathcal{J}|}\to \R$ for which the two expectations involving $\psi$ below are finite, we have
\[
\EE\left[\prod_{j\in \mathcal{J}} |X_j|^{-2 \nu_j} \times \psi((X_i)_{i\in \mathcal{J}^{\complement}})\right] \leq \EE\left[\prod_{j\in \mathcal{J}} |X_j|^{-2 \nu_j}\right] \times \EE\left[\psi((X_i)_{i\in \mathcal{J}^{\complement}})\right].
\]
\end{theorem}

Since the map
\[
(x_i)_{i\in \mathcal{J}^{\complement}}\mapsto \left(\sum_{i\in \mathcal{J}^{\complement}} \frac{\nu_i}{\sum_{i\in \mathcal{J}^{\complement}}\nu_i}|x_i|\right)^{2\sum_{i\in \mathcal{J}^{\complement}}\nu_i}
\]
is convex whenever $\nu_i\in [0,\infty)$ for all $i\in \mathcal{J}^{\complement}$ and $\sum_{i\in \mathcal{J}^{\complement}}\nu_i\in [1/2,\infty)$, the following corollary is immediate by combining Theorem~\ref{thm:upper.bound.convex} with the weighted arithmetic-mean-geometric-mean inequality.

\begin{corollary}\label{cor:upper.bound}
Let $\bb{X}\sim \mathcal{N}_d(\bb{0}_d, \Sigma)$ for some positive definite matrix $\Sigma$, and let $\emptyset\neq \mathcal{J}\subsetneq \{1,\ldots,d\}$ be given. Assume that $\nu_j\in [0,1/2)$ for all $j\in \mathcal{J}$, that $\nu_i\in [0,\infty)$ for all $i\in \mathcal{J}^{\complement}$, and that $\sum_{i\in \mathcal{J}^{\complement}}\nu_i\in [1/2,\infty)$. Then, we have
\[
\EE\left[\prod_{j\in \mathcal{J}} |X_j|^{-2 \nu_j} \times \prod_{i\in \mathcal{J}^{\complement}} |X_i|^{2 \nu_i}\right] \leq \EE\left[\prod_{j\in \mathcal{J}} |X_j|^{-2 \nu_j}\right] \times \EE\left[\left(\sum_{i\in \mathcal{J}^{\complement}} \frac{\nu_i}{\sum_{i\in \mathcal{J}^{\complement}}\nu_i}|X_i|\right)^{2\sum_{i\in \mathcal{J}^{\complement}}\nu_i}\right].
\]
\end{corollary}

Thus, in the range where the positive-side product can be dominated by a convex function, the product estimate goes in the upper-bound direction, rather than in the lower-bound direction of the base GPI. This is the multidimensional analogue of the direction in~\eqref{eq:Zhou.et.al.Thm.1.3} for the product cases covered by the corollary.

\begin{remark}\label{rem:2}
The argument leading to Corollary~\ref{cor:upper.bound} does not fully cover the range $\smash{\sum_{i\in \mathcal{J}^{\complement}}} \, \nu_i\in(0,1/2)$. Indeed, when $|\mathcal{J}^{\complement}| = 1$, this range is already covered by~\eqref{eq:Zhou.et.al.Thm.1.3}, but the remaining product case in which $|\mathcal{J}^{\complement}| \geq 2$ and $\sum_{i\in \mathcal{J}^{\complement}}\nu_i\in(0,1/2)$ is not covered by the present argument. One possible approach would be to exploit the Bernstein function integral representation, $|x|^{2\nu} = \nu \int_0^{\infty} (1 - e^{-s x^2}) s^{-\nu - 1} \rd s / \Gamma(1 - \nu)$, which holds for $x\in \R$ and $\nu\in (0,1/2)$; see, e.g., \citet[Eq.~(1)]{MR2978140}. However, we were unable to turn this representation into a successful argument. A treatment of this remaining case therefore lies beyond the scope of our paper.
\end{remark}

In the two-dimensional case, the following proposition gives an elliptical comparison inequality. Its Gaussian specialization \eqref{eq:Gaussian.specialization} is a functional strengthening of the two-dimensional case of Corollary~\ref{cor:upper.bound}, where one coordinate carries a negative-side exponent and the other carries a positive-side exponent.

\begin{proposition}\label{eq:2d.opposite.GPI.generalization}
Suppose that the two-dimensional random vector $(X,Y)$ has an elliptical distribution, meaning that
\[
(X,Y) \stackrel{\mathrm{law}}{=} R \, \Sigma^{1/2} \bb{U},
\]
where $R$ is a nonnegative random variable, $\Sigma^{1/2}$ is the symmetric square root of some positive definite matrix $\Sigma$, and $\bb{U}$ is a random vector uniformly distributed on the unit circle $\mathbb{S}^1$ that is independent of $R$. Also, let $f$ and $g$ be measurable functions on $[0,\infty)$ that are nonincreasing and nondecreasing on $(0,\infty)$, respectively, and for which the displayed expectations below are finite. Let $D = \mathrm{diag}(\sqrt{\Sigma_{11}},\sqrt{\Sigma_{22}})$, and let $(X_0,Y_0) \smash{\stackrel{\mathrm{law}}{=}} R D \bb{U}$, where $R$ and $\bb{U}$ are as above. Then
\begin{equation}\label{eq:elliptical.inequality}
\EE[f(|X|) g(|Y|)] \leq \EE[f(|X_0|) g(|Y_0|)].
\end{equation}
In particular, in the centered Gaussian case, that is, when $R\sim \chi_2$ and hence $R \, \bb{U}\sim \mathcal{N}_2(\bb{0}_2,I_2)$, so that $(X,Y)\sim \mathcal{N}_2(\bb{0}_2,\Sigma)$, we have
\begin{equation}\label{eq:Gaussian.specialization}
\EE[f(|X|) g(|Y|)] \leq \EE[f(|X|)] \EE[g(|Y|)].
\end{equation}
\end{proposition}

\begin{remark}
In \eqref{eq:Gaussian.specialization}, by taking $f(t)=t^{-2\nu_1}$ and $g(t)=t^{2\nu_2}$ for $t\in(0,\infty)$, with $\nu_1\in [0,1/2)$ and $\nu_2\in (0,\infty)$, and assigning arbitrary finite values at zero, we recover the opposite GPI of \citet[Theorem~2.1]{MR4471184}.
\end{remark}

\section{Proofs}\label{sec:proofs}

The set $\mathcal{J}$ indexes the negative-power factors, while the positive-power factors in the lower bound and the convex function in the upper bound depend on the coordinates indexed by $\mathcal{J}^{\complement}$. Both main proofs first use the Laplace-transform identity \eqref{eq:completely.monotonicity} for $|X_j|^{-2\nu_j}$, $j\in \mathcal{J}$, to rewrite the mixed expectation as an integral of expectations under quadratic Gaussian tilts. For each fixed mixture parameter $\bb{s}$, renormalizing the tilted density yields a centered Gaussian vector with covariance $(\Sigma^{-1} + T_{\bb{s}})^{-1}$. From this common representation, the proof of the lower bound applies the GPI to the tilted subvector indexed by $\mathcal{J}^{\complement}$ and then bounds its marginal variances from below by the corresponding Schur-complement variances. The proof of the upper bound instead uses $(\Sigma^{-1} + T_{\bb{s}})^{-1} \preceq \Sigma$, which compares the tilted and original vectors in the convex order. Indices with zero exponent can be omitted or handled by approximation, so the reduction separates the negative-power coordinates, which generate the tilt, from the coordinates to which the GPI or convex-order argument is applied.

\begin{proof}[\bf Proof of Theorem~\ref{thm:lower.bound}]
We first prove the claim when $\nu_j\in (0,1/2)$ for all $j\in \mathcal{J}$. Re-index the components of $\bb{X}$ so that $\mathcal{J} = \{1,\ldots,k\}$ with $1\leq k\leq d - 1$ and $\mathcal{J}^{\complement}=\{k+1,\ldots,d\}$. For every $x\in \R\setminus\{0\}$ and $\nu\in (0,1/2)$, the identity
\begin{equation}\label{eq:completely.monotonicity}
|x|^{-2\nu} = \frac{1}{\Gamma(\nu)} \int_0^{\infty} e^{-s x^2} s^{\nu - 1} \rd s
\end{equation}
holds; see, e.g., \citet[Eq.~(2)]{MR2978140}. This identity can be applied below since Gaussian coordinates are nonzero almost surely. Applying~\eqref{eq:completely.monotonicity} to $|X_1|^{-2\nu_1},\ldots,|X_k|^{-2\nu_k}$ and using Tonelli's theorem yields
\begin{equation}\label{eq:lower.bound.begin}
\begin{aligned}
&\EE\left[\prod_{j=1}^k |X_j|^{-2\nu_j} \prod_{i=k+1}^d |X_i|^{2\nu_i}\right] \\
&\qquad= \frac{1}{\prod_{j=1}^k \Gamma(\nu_j)} \int_{(0,\infty)^k} \EE\left[\exp(-\bb{X}^{\top} T_{\bb{s}} \bb{X}/2) \prod_{i=k+1}^d |X_i|^{2\nu_i}\right] \prod_{j=1}^k s_j^{\nu_j - 1} \rd \bb{s},
\end{aligned}
\end{equation}
where $\bb{s}=(s_1,\ldots,s_k)$ and $T_{\bb{s}} = \mathrm{diag}(2s_1,\ldots,2s_k,0,\ldots,0)$. For fixed $\bb{s}$, the integrand can be viewed as an expectation under an exponentially tilted centered Gaussian law. Using a renormalization argument, we have
\[
\begin{aligned}
&\EE\left[\exp(-\bb{X}^{\top} T_{\bb{s}} \bb{X}/2) \prod_{i=k+1}^d |X_i|^{2\nu_i}\right] \\
&\qquad= \int_{\R^d} \exp(-\bb{x}^{\top} T_{\bb{s}} \bb{x}/2) \, \frac{\exp(-\bb{x}^{\top} \Sigma^{-1} \bb{x}/2)}{(2\pi)^{d/2} |\Sigma|^{1/2}} \prod_{i=k+1}^d |x_i|^{2\nu_i} \rd \bb{x} \\
&\qquad= \frac{|(\Sigma^{-1} + T_{\bb{s}})^{-1}|^{1/2}}{|\Sigma|^{1/2}} \int_{\R^d} \frac{\exp\{-\bb{x}^{\top} (\Sigma^{-1} + T_{\bb{s}}) \, \bb{x}/2\}}{(2\pi)^{d/2} |(\Sigma^{-1} + T_{\bb{s}})^{-1}|^{1/2}} \prod_{i=k+1}^d |x_i|^{2\nu_i} \rd \bb{x} \\
&\qquad= \frac{|(\Sigma^{-1} + T_{\bb{s}})^{-1}|^{1/2}}{|\Sigma|^{1/2}} \, \EE\left[\prod_{i=k+1}^d |Y_i^{(\bb{s})}|^{2\nu_i}\right],
\end{aligned}
\]
where $\smash{\bb{Y}^{(\bb{s})}} = (Y_1^{(\bb{s})},\ldots,Y_d^{(\bb{s})})\sim\mathcal{N}_d(\bb{0}_d,(\Sigma^{-1} + T_{\bb{s}})^{-1})$. Then, using the assumption that $\mathbf{GPI}_{d - k}((\nu_i)_{i>k})$ holds, we obtain
\begin{equation}\label{eq:until}
\begin{aligned}
\EE\left[\exp(-\bb{X}^{\top} T_{\bb{s}} \bb{X}/2) \prod_{i=k+1}^d |X_i|^{2\nu_i}\right]
&\geq \frac{|(\Sigma^{-1} + T_{\bb{s}})^{-1}|^{1/2}}{|\Sigma|^{1/2}} \, \prod_{i=k+1}^d \EE\left[|Y_i^{(\bb{s})}|^{2\nu_i}\right] \\
&= \EE\left[\exp(-\bb{X}^{\top} T_{\bb{s}} \bb{X}/2)\right] \, \prod_{i=k+1}^d \EE\left[|Y_i^{(\bb{s})}|^{2\nu_i}\right].
\end{aligned}
\end{equation}
Substituting the last bound into~\eqref{eq:lower.bound.begin}, we obtain
\begin{equation}\label{eq:lower.bound.before.final.step}
\begin{aligned}
&\EE\left[\prod_{j=1}^k |X_j|^{-2\nu_j} \prod_{i=k+1}^d |X_i|^{2\nu_i}\right] \\
&\qquad\geq \frac{1}{\prod_{j=1}^k \Gamma(\nu_j)} \int_{(0,\infty)^k} \EE\left[\exp(-\bb{X}^{\top} T_{\bb{s}} \bb{X}/2)\right] \prod_{i=k+1}^d \EE\left[|Y_i^{(\bb{s})}|^{2\nu_i}\right] \prod_{j=1}^k s_j^{\nu_j - 1} \rd \bb{s}.
\end{aligned}
\end{equation}
It remains to find a lower bound for $\EE[|Y_i^{(\bb{s})}|^{2\nu_i}]$ that is uniform in $\bb{s}$. To this end, note that, for $i > k$,
\[
\begin{aligned}
\EE\left[|Y_i^{(\bb{s})}|^{2\nu_i}\right]
&= 2 \int_0^{\infty} u^{2\nu_i} \frac{1}{\sqrt{2\pi \, \Var(Y_i^{(\bb{s})})}} \exp\left\{-\frac{u^2}{2 \, \Var(Y_i^{(\bb{s})})}\right\} \rd u \\
&= \frac{2^{\nu_i} \{\Var(Y_i^{(\bb{s})})\}^{\nu_i}}{\sqrt{\pi}} \int_0^{\infty} t^{\nu_i-1/2} e^{-t} \rd t
= \frac{\Gamma(\nu_i + 1/2)}{\sqrt{\pi}} 2^{\nu_i} \{\Var(Y_i^{(\bb{s})})\}^{\nu_i},
\end{aligned}
\]
by making the change of variable $t = u^2/(2 \, \Var(Y_i^{(\bb{s})}))$. Next, write $\Sigma$ in block form:
\[
\Sigma
=
\begin{bmatrix}
A & B\\
B^{\top} & C
\end{bmatrix},
\]
where $A$ is $k\times k$, $B$ is $k\times(d - k)$ and $C$ is $(d - k)\times(d - k)$. Let $S = \mathrm{diag}(2 s_1,\ldots,2 s_k)$ so that $T_{\bb{s}}=\mathrm{diag}(S,0_{(d-k) \times (d-k)})$. The matrix inversion formula for block partitions in Lemma~\ref{lem:inverse.block.matrix.symmetric}~(i) gives
\[
\Sigma^{-1} + T_{\bb{s}}
=
\begin{bmatrix}
A^{-1} + S + A^{-1} B (\Sigma / A)^{-1} B^{\top} A^{-1} & - A^{-1} B (\Sigma / A)^{-1} \\[1mm]
- (\Sigma / A)^{-1} B^{\top} A^{-1} & (\Sigma / A)^{-1}
\end{bmatrix}
\equiv
\begin{bmatrix}
P & Q \\
Q^{\top} & R
\end{bmatrix}.
\]
Then Lemma~\ref{lem:inverse.block.matrix.symmetric}~(ii) yields
\[
\begin{aligned}
(\Sigma^{-1} + T_{\bb{s}})^{-1}
&=
\begin{bmatrix}
(A^{-1} + S)^{-1} & \star \\[1mm]
\star & (\Sigma/A) + B^{\top} A^{-1} (A^{-1} + S)^{-1} A^{-1} B
\end{bmatrix}.
\end{aligned}
\]
Let $\bb{e}_j$ denote the $j$th standard basis vector. Since $B^{\top} A^{-1} (A^{-1} + S)^{-1} A^{-1} B$ is nonnegative definite and $\Var(X_i) = \Sigma_{ii}$, we have, for $i>k$,
\[
\begin{aligned}
\Var(Y^{(\bb{s})}_i)
&= \bb{e}_{i-k}^{\top} \{(\Sigma/A) + B^{\top} A^{-1} (A^{-1} + S)^{-1} A^{-1} B\} \bb{e}_{i-k} \\[1mm]
&\geq \bb{e}_{i-k}^{\top} (\Sigma/A) \bb{e}_{i-k} \\
&= \frac{(\Sigma/A)_{i-k,i-k}}{\Sigma_{ii}} \Var(X_i),
\end{aligned}
\]
and therefore,
\[
\begin{aligned}
\EE\left[|Y_i^{(\bb{s})}|^{2\nu_i}\right]
&\geq \left\{\frac{(\Sigma/A)_{i-k,i-k}}{\Sigma_{ii}}\right\}^{\nu_i} \frac{\Gamma(\nu_i + 1/2)}{\sqrt{\pi}} 2^{\nu_i} \{\Var(X_i)\}^{\nu_i} \\
&= \left\{\frac{(\Sigma/A)_{i-k,i-k}}{\Sigma_{ii}}\right\}^{\nu_i} \EE[|X_i|^{2\nu_i}].
\end{aligned}
\]
Here, with a slight abuse of notation, we write $(\Sigma/\Sigma_{\mathcal{J}\mathcal{J}})_{ii}$ for $(\Sigma/A)_{i-k,i-k}$, with $i$ understood as its original index in $\mathcal{J}^{\complement}$. After substituting this lower bound into~\eqref{eq:lower.bound.before.final.step} and then using~\eqref{eq:completely.monotonicity} together with Tonelli's theorem, we get
\[
\EE\left[\prod_{j\in\mathcal{J}}|X_j|^{-2\nu_j} \times \prod_{i\in\mathcal{J}^{\complement}}|X_i|^{2\nu_i}\right]
\geq \EE\left[\prod_{j\in \mathcal{J}} |X_j|^{-2\nu_j}\right] \times \prod_{i\in\mathcal{J}^{\complement}} \left\{\frac{(\Sigma/\Sigma_{\mathcal{J}\mathcal{J}})_{ii}}{\Sigma_{ii}}\right\}^{\nu_i} \EE\left[|X_i|^{2\nu_i}\right],
\]
which is the claim.

If some of the exponents $\nu_j$, $j\in \mathcal{J}$, are equal to $0$, apply the preceding argument with those exponents replaced by numbers in $(0,1/2)$ that tend to $0$. The desired inequality follows from the dominated convergence theorem applied to both sides.
\end{proof}

\begin{proof}[\bf Proof of Theorem~\ref{thm:upper.bound.convex}]
If all the exponents $\nu_j$, $j\in \mathcal{J}$, are equal to $0$, then the result is immediate. If only some of them are equal to $0$, remove those indices from $\mathcal{J}$ and replace $\psi$ by the convex function on the enlarged complement that ignores the removed coordinates. Therefore, it is enough to prove the result when $\nu_j\in (0,1/2)$ for all $j\in \mathcal{J}$. We first prove the result when $\psi$ is nonnegative. Re-index the components of $\bb{X}$ so that $\mathcal{J} = \{1,\ldots,k\}$ with $1\leq k\leq d - 1$ and $\mathcal{J}^{\complement}=\{k+1,\ldots,d\}$. Following the same steps as in the proof of Theorem~\ref{thm:lower.bound} up until \eqref{eq:until}, but without applying $\mathbf{GPI}_{d - k}((\nu_i)_{i>k})$, we obtain
\[
\begin{aligned}
&\EE\left[\prod_{j=1}^k |X_j|^{-2\nu_j} \times \psi((X_i)_{i > k})\right] \\
&\qquad= \frac{1}{\prod_{j=1}^k \Gamma(\nu_j)} \int_{(0,\infty)^k} \EE\left[\exp(-\bb{X}^{\top} T_{\bb{s}} \bb{X}/2) \psi((X_i)_{i > k})\right] \prod_{j=1}^k s_j^{\nu_j - 1} \rd \bb{s},
\end{aligned}
\]
where $\smash{\bb{Y}^{(\bb{s})}} \sim\mathcal{N}_d(\bb{0}_d,(\Sigma^{-1} + T_{\bb{s}})^{-1})$ and
\[
\begin{aligned}
\EE\left[\exp(-\bb{X}^{\top} T_{\bb{s}} \bb{X}/2) \psi((X_i)_{i > k})\right]
&= \frac{|(\Sigma^{-1} + T_{\bb{s}})^{-1}|^{1/2}}{|\Sigma|^{1/2}} \, \EE\big[\psi((Y_i^{(\bb{s})})_{i > k})\big] \\
&= \EE\left[\exp(-\bb{X}^{\top} T_{\bb{s}} \bb{X}/2)\right] \, \EE\big[\psi((Y_i^{(\bb{s})})_{i > k})\big].
\end{aligned}
\]
Since $\Sigma^{-1} \preceq \Sigma^{-1} + T_{\bb{s}}$ in the Loewner order, we have $(\Sigma^{-1} + T_{\bb{s}})^{-1} \preceq \Sigma$. Thus, by Lemma~\ref{lem:ccx} applied to the convex function $\bb{x}\mapsto \psi((x_i)_{i > k})$, we have
\[
\EE\big[\psi((Y_i^{(\bb{s})})_{i > k})\big] \leq \EE\left[\psi((X_i)_{i > k})\right].
\]
Combining the last three displays and applying~\eqref{eq:completely.monotonicity} again together with Tonelli's theorem, we find that
\[
\EE\left[\prod_{j=1}^k |X_j|^{-2\nu_j} \times \psi((X_i)_{i > k})\right] \leq \EE\left[\prod_{j=1}^k |X_j|^{-2 \nu_j}\right] \times \EE\left[\psi((X_i)_{i > k})\right].
\]
This proves the result for nonnegative convex functions. For a general convex $\psi:\R^{d-k}\to \R$, let $\ell$ be an affine minorant of $\psi$ and let $\phi=\psi-\ell$. Then $\phi$ is nonnegative and convex. Since $\ell$ is affine and $\nu_j<1/2$ for all $j\leq k$, the corresponding expectations involving $\ell$ are finite, and hence the assumed finiteness of the two expectations involving $\psi$ implies that the two expectations involving $\phi$ are finite. Therefore, the preceding inequality holds with $\phi$ in place of $\psi$. For the affine part, the constant term is immediate, and the linear terms vanish by centeredness and invariance of both the Gaussian law and $\prod_{j=1}^k |X_j|^{-2\nu_j}$ under $\bb{x}\mapsto -\bb{x}$. Hence
\[
\EE\left[\prod_{j=1}^k |X_j|^{-2\nu_j} \times \ell((X_i)_{i > k})\right] = \EE\left[\prod_{j=1}^k |X_j|^{-2 \nu_j}\right] \times \EE\left[\ell((X_i)_{i > k})\right].
\]
Adding this identity to the inequality for $\phi$ proves the desired inequality for $\psi$. This concludes the proof.
\end{proof}

\begin{proof}[\bf Proof of Proposition~\ref{eq:2d.opposite.GPI.generalization}]
We first prove the result when $f$ and $g$ are bounded. Let $\bb{\alpha} = (\bb{e}_1^{\top} \Sigma^{1/2})^{\top}$ and $\bb{\beta} = (\bb{e}_2^{\top} \Sigma^{1/2})^{\top}$ denote, in column form, the first and second rows of $\Sigma^{1/2}$, respectively. Note that
\[
\|\bb{\alpha}\|_2 = \sqrt{(\bb{e}_1^{\top} \Sigma^{1/2}) (\bb{e}_1^{\top} \Sigma^{1/2})^{\top}} = \sqrt{\Sigma_{11}}, \qquad \|\bb{\beta}\|_2 = \sqrt{(\bb{e}_2^{\top} \Sigma^{1/2}) (\bb{e}_2^{\top} \Sigma^{1/2})^{\top}} = \sqrt{\Sigma_{22}}.
\]
Also, let $\bb{\alpha}_0 = \bb{\alpha}/\|\bb{\alpha}\|_2$ and $\bb{\beta}_0 = \bb{\beta}/\|\bb{\beta}\|_2$. When $r=0$, the two sides of the desired conditional inequality are both equal to $f(0)g(0)$. For fixed $r\in (0,\infty)$, define $F_r(x) = f(r \|\bb{\alpha}\|_2 \, x)$ and $G_r(y) = g(r \|\bb{\beta}\|_2 \, y)$. Then $F_r$ is nonincreasing and $G_r$ is nondecreasing on $(0,\infty)$. It is sufficient to prove
\[
\EE[F_r(|\bb{\alpha}_0 \meddot \bb{U}|) G_r(|\bb{\beta}_0 \meddot \bb{U}|)] \leq \EE[F_r(|\bb{e}_1 \meddot \bb{U}|) G_r(|\bb{e}_2 \meddot \bb{U}|)]
\]
for all $r\in (0,\infty)$, where $\bb{v}\meddot \bb{w} = v_1 w_1 + v_2 w_2$ denotes the dot product in $\R^2$.

To this end, let $T$ be the $2\times 2$ matrix representing a $\pi/2$ counterclockwise rotation, and define the function $H$ on $\mathbb{S}^1$ by
\[
H(\bb{u}) = \{F_r(|\bb{\alpha}_0 \meddot \bb{u}|) - F_r(|T \bb{\beta}_0 \meddot \bb{u}|)\} \{G_r(|\bb{\beta}_0 \meddot \bb{u}|) - G_r(|T \bb{\alpha}_0 \meddot \bb{u}|)\}.
\]
By Pythagoras, we have, for any $\bb{u}\in \mathbb{S}^1$,
\[
|\bb{\alpha}_0 \meddot \bb{u}|^2 + |T \bb{\alpha}_0 \meddot \bb{u}|^2 = |\bb{\beta}_0 \meddot \bb{u}|^2 + |T \bb{\beta}_0 \meddot \bb{u}|^2 = \|\bb{u}\|_2^2 = 1.
\]
It follows that
\[
|\bb{\alpha}_0 \meddot \bb{u}|^2 - |T \bb{\beta}_0 \meddot \bb{u}|^2 = |\bb{\beta}_0 \meddot \bb{u}|^2 - |T \bb{\alpha}_0 \meddot \bb{u}|^2,
\]
which in turn implies that $|\bb{\alpha}_0 \meddot \bb{u}| \geq |T \bb{\beta}_0 \meddot \bb{u}| ~\Leftrightarrow~ |\bb{\beta}_0 \meddot \bb{u}| \geq |T \bb{\alpha}_0 \meddot \bb{u}|$. Since the four inner products in the definition of $H$ vanish only on a finite subset of $\mathbb{S}^1$, and since $F_r$ is nonincreasing and $G_r$ is nondecreasing on $(0,\infty)$, we have $H(\bb{u}) \leq 0$ for almost every $\bb{u}\in \mathbb{S}^1$, and thus
\[
\EE[H(\bb{U})] \leq 0.
\]
On the other hand, since the ordered pairs $(\bb{\alpha}_0,\bb{\beta}_0)$ and $(T \bb{\beta}_0,T \bb{\alpha}_0)$ have the same Gram matrix, the orthogonal invariance of $\bb{U}$ gives
\[
\EE[F_r(|\bb{\alpha}_0 \meddot \bb{U}|) G_r(|\bb{\beta}_0 \meddot \bb{U}|)] = \EE[F_r(|T \bb{\beta}_0 \meddot \bb{U}|) G_r(|T \bb{\alpha}_0 \meddot \bb{U}|)].
\]
Similarly, since $\bb{\alpha}_0 \meddot T \bb{\alpha}_0 = T \bb{\beta}_0 \meddot \bb{\beta}_0 = \bb{e}_1 \meddot \bb{e}_2 = 0$, we have
\[
\begin{aligned}
\EE[F_r(|\bb{\alpha}_0 \meddot \bb{U}|) G_r(|T \bb{\alpha}_0 \meddot \bb{U}|)]
&= \EE[F_r(|T \bb{\beta}_0 \meddot \bb{U}|) G_r(|\bb{\beta}_0 \meddot \bb{U}|)] \\
&= \EE[F_r(|\bb{e}_1 \meddot \bb{U}|) G_r(|\bb{e}_2 \meddot \bb{U}|)].
\end{aligned}
\]
Combining the last three displays yields
\[
\EE[F_r(|\bb{\alpha}_0 \meddot \bb{U}|) G_r(|\bb{\beta}_0 \meddot \bb{U}|)] \leq \EE[F_r(|\bb{e}_1 \meddot \bb{U}|) G_r(|\bb{e}_2 \meddot \bb{U}|)].
\]
Equivalently,
\[
\EE[f(r|\bb{\alpha}\meddot \bb{U}|) g(r|\bb{\beta}\meddot \bb{U}|)] \leq \EE[f(r \|\bb{\alpha}\|_2 \, |U_1|) g(r \|\bb{\beta}\|_2 \, |U_2|)].
\]
Together with the equality at $r=0$, this inequality holds for every $r\in[0,\infty)$. Integrating it with respect to the distribution of $R$ gives
\[
\EE[f(|X|) g(|Y|)] \leq \EE[f(|X_0|) g(|Y_0|)].
\]
The general elliptical inequality \eqref{eq:elliptical.inequality} follows by applying the preceding bounded case to the truncations
\[
f_n(t)=\max\{-n,\min\{f(t),n\}\}, \qquad g_n(t)=\max\{-n,\min\{g(t),n\}\},
\]
which preserve the required monotonicity on $(0,\infty)$, and then using dominated convergence, since $|f_n(x)g_n(y)|\leq |f(x)g(y)|$ for all $x,y\in[0,\infty)$.

Finally, assume that $R\sim \chi_2$. Then $R \, \bb{U}\sim \mathcal{N}_2(\bb{0}_2,I_2)$, so $(X,Y)\sim \mathcal{N}_2(\bb{0}_2,\Sigma)$. Moreover, $(X_0,Y_0)\sim \mathcal{N}_2(\bb{0}_2,D^2)$, and therefore $X_0$ and $Y_0$ are independent. Since $X_0 \smash{\stackrel{\mathrm{law}}{=}} X$ and $Y_0 \smash{\stackrel{\mathrm{law}}{=}} Y$, it follows that
\[
\EE[f(|X_0|) g(|Y_0|)] = \EE[f(|X_0|)] \EE[g(|Y_0|)] = \EE[f(|X|)] \EE[g(|Y|)].
\]
Together with \eqref{eq:elliptical.inequality}, this proves the Gaussian specialization \eqref{eq:Gaussian.specialization}.
\end{proof}

\begin{appendices}

\renewcommand{\thesection}{Appendix~\Alph{section}}

\section{Technical lemmas}\label{app:tech.lemmas}

\renewcommand{\thesection}{\Alph{section}}

The first lemma contains well-known formulas from matrix analysis for the inverse of a $2\times 2$ block matrix; see, e.g., Theorem~2.1 of \citet{MR1873248}.

\begin{lemma}\label{lem:inverse.block.matrix.symmetric}
Let $\Sigma$ and $M$ be symmetric matrices of width at least $2$, each written in $2\times 2$ block form.
\begin{itemize}
\item[(i)]
If $A$ is invertible and the Schur complement $\Sigma / A = C - B^{\top} A^{-1} B$ is invertible, then
\[
\Sigma =
\begin{bmatrix}
A & B \\
B^{\top} & C
\end{bmatrix}
\quad \Rightarrow \quad
\Sigma^{-1} =
\begin{bmatrix}
A^{-1} + A^{-1} B (\Sigma / A)^{-1} B^{\top} A^{-1} & - A^{-1} B (\Sigma / A)^{-1} \\[1mm]
- (\Sigma / A)^{-1} B^{\top} A^{-1} & (\Sigma / A)^{-1}
\end{bmatrix}.
\]
\item[(ii)]
If $R$ is invertible and the Schur complement $M / R = P - Q R^{-1} Q^{\top}$ is invertible, then
\[
M =
\begin{bmatrix}
P & Q \\
Q^{\top} & R
\end{bmatrix}
\quad \Rightarrow \quad
M^{-1} =
\begin{bmatrix}
(M / R)^{-1} & - (M / R)^{-1} Q R^{-1} \\[1mm]
- R^{-1} Q^{\top} (M / R)^{-1} & R^{-1} + R^{-1} Q^{\top} (M / R)^{-1} Q R^{-1}
\end{bmatrix}.
\]
\end{itemize}
\end{lemma}

The second lemma states that multivariate normal random vectors are ordered in the convex order according to the Loewner order of their covariance matrices; see, e.g., \citet[Example~7.A.13]{MR2265633}.

\begin{lemma}\label{lem:ccx}
Let $\Sigma_1,\Sigma_2$ be $d\times d$ positive definite matrices, and suppose that $\bb{Y}\sim \mathcal{N}_d(\bb{0}_d,\Sigma_1)$ and $\bb{X}\sim \mathcal{N}_d(\bb{0}_d,\Sigma_2)$. If $\Sigma_1 \preceq \Sigma_2$, then $\bb{Y}$ is less than or equal to $\bb{X}$ in the convex order, written $\bb{Y}\leq_{\mathrm{cx}} \bb{X}$, meaning that for every convex function $\psi:\R^d\to \R$ for which the expectations below are finite, we have
\[
\EE[\psi(\bb{Y})] \leq \EE[\psi(\bb{X})].
\]
\end{lemma}

\end{appendices}

\section*{Funding}
\addcontentsline{toc}{section}{Funding}

Guolie Lan's research is supported by the National Natural Science Foundation of China (grant No.\ 12171335) and the Guangdong Provincial Natural Science Foundation (grant No.\ 2023A1515012170). Fr\'ed\'eric Ouimet acknowledges funding from the Natural Sciences and Engineering Research Council of Canada (NSERC) through Discovery grant RGPIN-2026-04471 and Discovery Launch Supplement DGECR-2026-00449. An early version of this work was completed while Fr\'ed\'eric Ouimet was a Research Associate at McGill University and a postdoctoral fellow at the Université de Sherbrooke. These positions were funded through Christian Genest's research grants (NSERC grant RGPIN-2024-04088 and Canada Research Chairs Program grant 950-231937) and Anne MacKay's NSERC grant (RGPIN-2024-05794). Wei Sun's research is funded by NSERC (grant RGPIN-2025-06779).

\section*{Acknowledgments}
\addcontentsline{toc}{section}{Acknowledgments}

We are grateful to the referees for their careful reading of the manuscript and for their constructive comments, which helped us improve the presentation and clarity of the paper.

\addcontentsline{toc}{section}{References}


\end{document}